\documentclass[twoside,11pt]{article}
\title{} \author{} \date{}
\usepackage{latexsym,amssymb,times,graphicx}%,showkeys}
\usepackage[table]{xcolor}% http://ctan.org/pkg/xcolor
\input amssym.def
\newtheorem{te}{Theorem}[section]
\newtheorem{prop}[te]{Proposition}

\newtheorem{fac}[te]{Fact}

\newtheorem{lem}[te]{Lemma}

\newtheorem{rem}[te]{Remark}
\newtheorem{ex}[te]{Example}

\def\dok{\noindent{\bf Proof. }}
\def\kdok{\hfill $\Box$ \par \vspace*{2mm} }

\let\strokel\l

\def\a{\alpha}

\def\g{\gamma}
\def\d{\delta}

\def\o{\omega}
\def\k{\kappa}
\def\l{\lambda}

\def\r{\rho}
\def\s{\sigma}

\def\t{\tau}

\def\S{{\mathbb S}}

\def\P{{\mathbb P}}
\def\Q{{\mathbb Q}}

\def\N{{\mathbb N}}
\def\X{{\mathbb X}}
\def\Y{{\mathbb Y}}
\def\Z{{\mathbb Z}}

\def\A{{\mathbb A}}

\def\BL{{\mathbb L}}

\def\L{{\mathcal L}}

\def\CX{{\mathcal X}}

\def\CZ{{\mathcal Z}}
\def\CW{{\mathcal W}}

\def\la{\langle}
\def\ra{\rangle}

\def\otp{\mathop{\rm otp}\nolimits}

\def\Scatt{\mathop{\mbox{\rm Scatt}}\nolimits}

\def\Sym{\mathop{\rm Sym}\nolimits}

\def\Lim{\mathop{\rm Lim}\nolimits}
\def\Sur{\mathop{\rm Sur}\nolimits}
\def\rfm{\mathop{\rm RFM}\nolimits}

\def\Ord{\mathop{\rm Ord}\nolimits}

\begin{document}

\thispagestyle{plain}
\begin{center}
           {\large \bf \uppercase{Reversible disjoint unions of well orders \\[1mm]and their inverses}}
\end{center}
\begin{center}
{\bf Milo\v s S.\ Kurili\'c\footnote{Department of Mathematics and Informatics, Faculty of Science, University of Novi Sad,
              Trg Dositeja Obradovi\'ca 4, 21000 Novi Sad, Serbia.
              email: milos@dmi.uns.ac.rs}
and Nenad Mora\v ca\footnote{Department of Mathematics and Informatics, Faculty of Science, University of Novi Sad,
              Trg Dositeja Obradovi\'ca 4, 21000 Novi Sad, Serbia.
              email:
              nenad.moraca@dmi.uns.ac.rs}}
\end{center}
\begin{abstract}
\noindent
A poset $\P$ is called reversible iff every bijective homomorphism $f:\P \rightarrow \P$ is an automorphism.
Let $\CW$ and $\CW ^*$ denote the classes of well orders and their inverses respectively.
We characterize reversibility in the class of posets of the form $\P =\bigcup _{i\in I}\BL _i$, where $\BL _i, i\in I$, are pairwise disjoint
linear orders from $\CW \cup \CW ^*$. First, if $\BL _i \in \CW$, for all $i\in I$, and $\BL _i \cong \a _i =\gamma_i+n_i\in \Ord$, where
$\gamma_i\in \Lim \cup \{0\}$ and $n_i\in\omega$, defining
$I_\alpha   :=  \{i\in I:\alpha_i=\alpha\}$, for $\a\in \Ord$, and
$J_\gamma   :=  \{j\in I:\gamma_j=\gamma\}$, for $\g\in \Lim _0$,
we prove that
$\bigcup _{i\in I} \BL_i$ is a reversible poset iff
 $\la \a_i :i\in I\ra$ is a finite-to-one sequence,
 or there is $\g =\max \{ \g _i : i\in I\}$, for $\a \leq \g$ we have $|I_\a|<\o  $,
  and
  $\la n_i : i\in J_\g \setminus I_\g\ra$ is a reversible sequence of natural numbers.
The same holds when $\BL _i \in \CW^*$, for all $i\in I$. In the general case, the reversibility of the whole union
is equivalent to the reversibility of the union of components from $\CW$ and the union of components from $\CW^*$.


\end{abstract}
\section{Introduction}\label{S1}
A relational structure $\X=\la X,\r\ra $ is called {\it reversible}
iff $ \la X,\s\ra\ \not\cong \la X,\r\ra $, whenever $\s\varsubsetneq \r$
(that is, iff every bijective endomorphism $f:\X\rightarrow \X$ is an automorphism of $\X$).
The class of reversible structures contains several relevant classes of structures; for example,
linear orders and, more generally, all structures first-order definable in linear orders by quantifier-free formulas without parameters
(i.e., monomorphic or chainable structures) \cite{KDef}, tournaments, Henson graphs \cite{KuMoExtr}, and Henson digraphs \cite{KRet}.
In addition, reversibility is preserved under some forms of bi-interpretability \cite{KRet} (but  bi-definability,
bi-embedability and elementary equivalence  do not preserve it in general \cite{KDef,KuMoSim})
and extreme elements of $L_{\infty \o}$-definable classes of structures are reversible under certain syntactical restrictions \cite{KuMoExtr}.
We also note that reversible structures have the Cantor-Schr\"{o}der-Bernstein property for condensations
(bijective homomorphisms). Namely, writing $\mathbb{X}\preccurlyeq_c\mathbb{Y}$ iff there is a condensation
$f:\mathbb{X}\rightarrow\mathbb{Y}$, we have: if $\X$ is a reversible structure, $\X \preccurlyeq_c \Y$ and $\Y \preccurlyeq_c \X$, then $\Y\cong \X$.
%and there are condensations $f:\X \rightarrow \Y$ and $g:\Y \rightarrow \X$, then there is an
%isomorphism $F:\X \rightarrow \Y$.

Of course, some (elementary) classes of structures contain both reversible and non-reversible structures; an example is
the class of equivalence relations, where reversibility was characterized in \cite{KuMoEqu}.
It seems that the corresponding characterization for the class of partial orders is a much more difficult task.
In \cite{Kuk} Kukie\strokel a initiated the investigation on that topic detecting some subclasses of reversible posets;
for example, Boolean lattices and well founded posets with finite levels, and in \cite{Kuk1} characterized
hereditarily reversible posets (the posets having all substructures reversible).
The article \cite{KuMoDiscI} contains several equivalents of reversibility in the class of disconnected binary structures.
First, it is easy to see that the connectivity components of a reversible structure have to be reversible; more generally we have
\begin{fac}\label{TC003}
If $\mathbb{X}_i, i\in I$ are pairwise disjoint and connected  binary structures, then
$\bigcup _{i\in I} \mathbb{X}_i$ is reversible iff $\;\bigcup_{i\in J}\mathbb{X}_i$ is reversible for each non-empty set $J\subset I$.
\end{fac}
Second, roughly speaking and denoting the order type of the integer line $\Z$ by $\zeta$, a structure having reversible components is reversible iff
its components can not be ``merged" by condensations and each $\zeta$-sequence of condensations between different components must be a sequence of isomorphisms.
A consequence of that characterization is the following statement from \cite{KuMoDiscI},
which will be used in this article. Denoting by $\Sym (I)$ (resp.\ $\Sur (I)$) the set of all bijections (resp.\ surjections) $f:I\rightarrow I$,
we will call a sequence of pairwise disjoint binary structures $\langle\mathbb{X}_i:i\in I\rangle$  a {\it reversible sequence of structures} iff
\begin{equation}\label{EQC007}\textstyle
\neg\exists f\in\Sur(I)\setminus\Sym(I)\;\;\forall j\in I\;\;\bigcup_{i\in f^{-1}[\{j\}]}\mathbb{X}_i\preccurlyeq_c\mathbb{X}_j.
\end{equation}
\begin{fac}\label{TC002}
If $\,\mathbb{L}_i$, $i\in I$, are pairwise disjoint linear orders, then  the following conditions are equivalent:

(a) The poset $\bigcup_{i\in I}\mathbb{L}_i$ is reversible;

(b) $\langle\mathbb{L}_i:i\in I\rangle$ is a reversible sequence of structures;

(c) There is no $f\in\Sur(I)\setminus\Sym(I)$ such that each component $\mathbb{L}_j$ can be partitioned into copies of $\mathbb{L}_i$,
where $i\in f^{-1}[\{j\}]$.
\end{fac}
Third, several sufficient conditions for the reversibility of the posets of the form $\bigcup _{i\in I}\mathbb{L}_i$ are given in \cite{KuMoDiscI}.
For example, denoting the bi-embedability relation by $\rightleftarrows$ and
defining a linear order $\BL$ to be {\it Cantor-Schr\"{o}der-Bernstein (for embeddings)}, shortly CSB,  iff for each linear order $\BL'$
satisfying $\BL'\rightleftarrows \BL$ we have $\BL'\cong \BL$ (that is, iff $[\BL ]_{\rightleftarrows }=[\BL ]_{\cong}=:\otp (\BL )$) we have
\begin{fac}\label{TB055}
If $\BL _i, i\in I$, are pairwise disjoint $\s$-scattered linear orders,\footnote{A linear order $\BL$ is called
{\it scattered} iff $\Q \not\hookrightarrow \BL$, where $\Q$ is the rational line.
$\BL$ is said to be {\it $\s$-scattered} iff $\BL$ is at most countable union of scattered linear orders.}
then the poset $\bigcup _{i\in I}\BL_i$ is reversible if
\begin{equation}\label{EQC006}
\la [\BL _i]_\rightleftarrows :i\in I\ra \mbox{ is a finite-to-one sequence}.\footnote{that is,
there is no infinite $J\subset I$ such that $\BL _i \rightleftarrows \BL _j$, for all $i,j\in J$.}
\end{equation}
If, in addition $\BL _i$'s are CSB linear orders, then the poset $\bigcup _{i\in I}\BL_i$ is reversible if
\begin{equation}\label{EQC006'}
\la \otp (\BL _i ) :i\in I\ra \mbox{ is a finite-to-one sequence}.\footnote{that is,
there is no infinite $J\subset I$ such that $\BL _i \cong \BL _j$, for all $i,j\in J$.}
\end{equation}
\end{fac}
We recall that a scattered linear order $\BL$ is said to be {\it of a limit type}  iff $\BL \rightleftarrows \sum _{s\in S} \BL _s$, where
$\S\in \Scatt$ and
$\BL _s\cong \o$ or $\BL _s\cong\o^*$, for each $s\in S$. For characterizations of such chains see \cite{Lav1}, p.\ 112.
E.g., the linear order $\BL= \o\o^* +1$ is of a limit type because $\BL\rightleftarrows \o\o^*$ but, since $\BL\not\cong \o\o^*$, it is not
CSB. Successor ordinals are CSB,  but not of a limit type. Limit ordinals are CSB of a limit type.

Let
$\CW$ denote the class of well orders,
$\L $ the class of well orders isomorphic to limit ordinals,
$\CZ$ the class of linear orders isomorphic to $\o ^\theta \o ^* +\o ^{\d }$, where $\theta$ and $\d$ are ordinals satisfying $1\leq\theta<\d$,
and let
$\CW ^*$, $\L ^*$, and $\CZ ^*$ be the classes of the inverses of elements of $\CW $, $\L $, and $\CZ $, respectively.
In \cite{Laf} Laflamme, Pouzet and Woodrow have shown that a scattered linear order is CSB iff it is isomorphic to a finite sum of linear orders from
$\CW \cup \CW ^* \cup \CZ \cup \CZ ^*$.
Also, by Theorem 5.4 of \cite{KuMoDiscI},
\begin{fac}\label{TC008}
(a) A linear order is CSB of a limit type iff it is isomorphic to a finite sum of linear orders from
$\L \cup \L ^* \cup \CZ \cup \CZ ^*$.

(b) If $\BL _i, i\in I$, are pairwise disjoint CSB linear orders of a limit type, then the poset $\bigcup _{i\in I}\BL_i$
is reversible iff condition (\ref{EQC006'}) is satisfied.
\end{fac}
Thus, in particular, condition (\ref{EQC006'}) is equivalent to reversibility
in the class of posets of the form $\bigcup _{i\in I}\BL_i$, where $\BL _i\in \L \cup \L ^*$, for all $i\in I$,
and in this article we characterize reversibility
in the class of posets of the form $\bigcup _{i\in I}\BL_i$, where $\BL _i\in \CW \cup \CW ^*$, for all $i\in I$.
In Section \ref{S2} we consider the cases when $\BL _i\in \CW $, for all $i\in I$, or $\BL _i\in \CW ^*$, for all $i\in I$,
and in Section \ref{S3} we observe the general case.

The rest of this section contains additional definitions and basic facts which will be used in the paper.
A sequence $\la n _i :i\in I\ra$ in $\N$ is defined to be a
{\it reversible sequence of natural numbers} iff
\begin{equation}\label{EQB032}\textstyle
\neg \exists f\in \Sur (I)\setminus \Sym (I)\;\;\forall j\in I \;\;\sum _{i\in f^{-1}[\{ j \}]}n_i =n _j .
\end{equation}
In order to give a characterization of such sequences from \cite{KuMoEqu} we recall that a set $K\subset\N$ %of natural numbers
is called {\it independent}
iff $n\not\in \la K\setminus \{ n \}\ra$, for all $n\in K$,
where $\la K\setminus \{ n \}\ra$ is the subsemigroup of the semigroup  $\la \N , +\ra$
generated by $K\setminus \{ n \}$; thus, $\emptyset$ is an independent set. By $\gcd (K)$  we denote the greatest common divisor of the numbers from $K$.
Defining $I_m:=\{ i\in I : n_i=m\}$, for $m\in \N$, by \cite{KuMoEqu} we have
\begin{fac}\label{TC001}
A sequence $\langle n_i :i\in I\rangle \in {}^{I}\N$ is reversible iff
%$K=\{ m\in \mathbb{N} : n_i =m ,\mbox{ for infinitely many } i\in I \}$
$K=\{ m\in \mathbb{N} : |I_m|\geq \o\}$
is an independent set and, if $K\neq \emptyset$, then $\gcd (K)$ divides at most finitely many elements
of the set $\{ n_i :i\in I \}$.\footnote{For example, if $I$ is a non-empty set of any size
and $\la n_i :i\in I\ra \in {}^{I}\N$, then by Fact \ref{TC001} we have:
 if $K=\emptyset$ (which is possible if $|I|\leq \o$), then $\la n_i \ra $ is a reversible sequence;
 if $K=\{ 2,5\}$, then  $\la n_i \ra $ is a reversible sequence iff the set $\{ n_i :i\in I \}$ is finite;
 if $K=\{ 4,10\}$, then  $\la n_i \ra $ is a reversible sequence iff the set $\{ n_i :i\in I \}$
contains at most finitely many even numbers.
}
\end{fac}
Let $\Ord$ and $\Lim$ denote the classes of ordinals and limit ordinals respectively.
\begin{fac}\label{TC000}
For each $\g\in \Lim$ we have

(a) If $1\leq\l\leq \o$, then there is a partition $\{ A_k :k< \l\}$ of $\g$
such that $A_k \cong \g $, for all $k< \l$;

(b) If $m,n\in \o$ and $f:\g +m \hookrightarrow \g +n$, then $f[\g ]$ is an unbounded subset of $\g$ and $f[m]\subset n$;

(c) If $\a \leq \g$, there is a partition $\{ A,B\}$ of $\g$, where $A\cong \a$ and $B\cong \g$.
\end{fac}
\dok
(a) We recall (see \cite{Rosen}, p.\ 71) that if $\BL =\la L, <\ra$ is a linear order and $\sim$  the equivalence relation on $L$ defined by
$x\sim y $ iff $|[\min\{x,y\},\max\{x,y\} ]|<\o$,  then  $\BL =\sum _{t\in T}\BL _t$, where  $\{ L_t :t\in T\}=L/\!\sim $ and
$\otp (\BL _t)\in \N \cup \{ \o ,\o ^* ,\zeta\}$, for each $t\in T$.
So we have $\g =\sum _{t\in T}\BL _t$, where $\BL_t\cong \o$.
For $t\in T$ we choose a partition $\{ A^k_t :k<\l\}$ of $L_t$ such that $|A^k_t|=\o$, for all $k<\l$.
Then, clearly, $A^k_t \cong \o$ and, hence, $A_k:=\sum _{t\in T}A^k_t \cong \sum _{t\in T}\BL _t=\g$.

(b) Since $\g\in \Lim$ we have $f[\g ]\subset \g$. Assuming that $f[\g ]\subset \a$, for some $\a<\g$, we would have
$\g =\otp f[\g ]\leq \a<\g$, which is false. Thus $f[\g ]$ is unbounded in $\g$
and, hence, $f[m]\subset n$.

(c) By (a), there is a partition $\{ A_0 ,A_1\}$ of $\g$ such that $A_0\cong A_1 \cong \g $. Since $\a \leq \g$
there is $A\subset A_0$ such that $A\cong \a$ and, clearly, $\g\cong A_1 \subset B:= \g \setminus A \subset \g$,
which implies that $B\cong \g$.
\kdok
\section{Reversible unions of well orders}\label{S2}
In the sequel we assume that $I$ is a non-empty set and $\langle\alpha_i:i\in I\rangle$ an $I$-sequence of nonzero ordinals. For $i\in I$,
$\alpha_i=\gamma_i+n_i$ will be the unique decomposition of $\alpha_i$, where $\gamma_i\in \Lim _0 :=\Lim \cup \{ 0\}$ and $n_i\in\omega$.

In addition, we assume that $\BL _i, i\in I$, are pairwise disjoint well orders, where $\BL _i=\BL _i' +\BL _i''$, $\BL _i'\cong \g _i$, $\BL _i''\cong n _i$, and, hence,
$\BL _i\cong \a _i$ and we consider the poset $\P =\bigcup _{i\in I}\BL _i $, which will be denoted by $\bigcup _{i\in I}\a _i$ or $\bigcup _{i\in I}\g _i +n_i$,
whenever there is no danger of confusion.

Also, we define
\begin{eqnarray*}
I_\alpha  & := & \{i\in I:\alpha_i=\alpha\}, \mbox{ for } \a\in \Ord ,\\
J_\gamma  & := & \{j\in I:\gamma_j=\gamma\}, \mbox{ for } \g\in \Lim _0,
\end{eqnarray*}
and in the sequel prove the following characterization.
\begin{te}\label{TC006}
$\bigcup _{i\in I} \a_i$ is a reversible poset iff exactly one of the following is true
\begin{itemize}
\item[\rm(I)] $\la \a_i :i\in I\ra$ is a finite-to-one sequence,
%\item[(D)] There is $\g =\max \{ \g _i : i\in I\}$,  $|I_\a |<\o$, for all $\a \leq \g$, $|J_\g|\geq \o$ and
%           $\la n_i : i\in J_\g \setminus I_\g\ra$ is a reversible sequence, which is not finite-to-one.
%
\item[\rm(II)] There is $\g =\max \{ \g _i : i\in I\}$, \\
           for $\a \leq \g$ we have $|I_\a|<\o  $,\\ %, $|J_\g|\geq \o$
           and the sequence $\la n_i : i\in J_\g \setminus I_\g\ra$ is reversible, but not finite-to-one.
\end{itemize}
\noindent
The same holds for the poset  $\bigcup _{i\in I} \a_i ^*$.
\end{te}
\begin{lem}\label{TC004}
The poset $\bigcup_{i\in I}\g +n_i$ is reversible iff $|I_\g|< \o$ and  $\langle n_i:i\in I\setminus I_\g\rangle$
is a reversible sequence of natural numbers.
\end{lem}
\dok
Let $\P =\bigcup_{i\in I}\BL _i$, where $\BL _i =\BL _i'+ \BL _i''\cong \g +n_i$, for $i\in I$, are disjoint linear orders.
We prove the contrapositive of the statement.

Let the poset $\P$ be non-reversible and $|I_\g|<\o$. Then by Fact \ref{TC002} there is $f\in \Sur (I)\setminus \Sym (I)$
such that for each $j\in J$ there is a
partition $\{ A_i: i\in f^{-1}[\{ j\}]\}$ of $L_j$
such that $\A_i = \A_i' +\A_i''$, where $\A_i'\cong \mathbb{L} _i'$ and $\A_i''\cong \mathbb{L} _i''$, for all $i\in f^{-1}[\{ j\}]$.

If $i\in I\setminus I_\g$, then $\BL _i$ is a successor ordinal bigger than $\g$
and, since $\BL _i\hookrightarrow \BL _{f(i)}$, we have $f(i)\in I\setminus I_\g$. So $f[ I\setminus I_\g]\subset I\setminus I_\g$
and, since $|f[I_\g]|<\o$ and $f$ is a surjection, $f[ I\setminus I_\g]= I\setminus I_\g$, which implies that
\begin{equation}\label{EQC000}
f[ I_\g]= I_\g \;\mbox{ and }\; f\upharpoonright (I\setminus I_\g)\in \Sur (I\setminus I_\g)\setminus \Sym (I\setminus I_\g).
\end{equation}
If $j\in I\setminus I_\g$, $i\in f^{-1}[\{ j\}]$ and $g_i:\BL _i \hookrightarrow \BL _j$, where $g_i[L_i]=A_i$,
then by (\ref{EQC000}) we have $i\in I\setminus I_\g$ and, by Fact \ref{TC000}(b),  $g_i[L_i']=A _i'\subset L_j'$ and $g_i[L_i'']=A _i''\subset L_j''$.
Since  $\{ A_i: i\in f^{-1}[\{ j\}]\}$ is a partition of $L_j$
it follows that $\{ A_i'': i\in f^{-1}[\{ j\}]\}$ is a partition of $L_j''$.
Now, since $\A_i''\cong n_i$ and $\BL_j''\cong n_j$
we have $\sum _{i\in f^{-1}[\{ j \}]}n_i =n _j$, for all $j\in I\setminus I_\g$, and $\langle n_i:i\in I\setminus I_\g\rangle$
is not a reversible sequence.

Conversely, if $|I_\g|\geq\o$, then $\g >0$ and, by Fact \ref{TC008}(b) $\bigcup _{i\in I _\g}\BL _i$ is not reversible;
so $\P$ is not reversible by Fact \ref{TC003}.

If $\langle n_i:i\in I\setminus I_\g\rangle$ is not a reversible sequence in $\N$ we show that the union
$\bigcup_{i\in I\setminus I_\g}\BL_i$ is not a reversible structure, which will, by Fact \ref{TC003}, show that $\P$ is not reversible as well.
W.l.o.g.\ assume that $I_\g=\emptyset$
and let $f\in \Sur (I)\setminus \Sym (I)$, where
$\sum _{i\in f^{-1}[\{ j \}]}n_i =n _j$, for all $j\in I $. Then for $j\in J$ we have $|f^{-1}[\{ j \}]|<\o$
and by Fact \ref{TC000}(a) there is a partition $\{ A_i' : i\in f^{-1}[\{ j \}]\}$ of the set $L_{j}'\cong \g$
such that $A_i'\cong \g\cong \BL _i'$, for all $i\in f^{-1}[\{ j \}]$.
Since $\sum _{i\in f^{-1}[\{ j \}]}n_i =n _j$ there is a partition $\{ A_i'' : i\in f^{-1}[\{ j \}]\}$ of the set $L_{j}''\cong n_j$
such that $A_i''\cong n_i\cong \BL _i''$, for all $i\in f^{-1}[\{ j \}]$.
Now defining $A_i:=A_i'\cup A_i''$, for $i\in f^{-1}[\{ j \}]$, we obtain a partition
$\{ A_i : i\in f^{-1}[\{ j \}]\}$ of $L_{j}$, where $A_i\cong \g +n_i\cong \BL _i$, for all $i\in f^{-1}[\{ j \}]$.
By Fact \ref{TC002} the structure $\bigcup _{i\in I}\mathbb{L}_i$ is not reversible.
\kdok
\begin{te}\label{TC005}
$\bigcup _{i\in I} \a _i$ is not reversible iff at least one of the following holds:
\begin{itemize}
\item[\rm (A)] There are $i\in I$ and $\a \leq \g _{i}$ such that $|I_{\a }|\geq \o $,
\item[\rm (B)] There is $\g \in \Lim_0$  such that $|J_\g \setminus I_\g|\geq \o$\\ and the sequence $\la n_i : i\in J_\g \setminus I_\g\ra$ is not reversible.
\end{itemize}
\end{te}
\dok
($\Leftarrow$) Let $i_0\in I$, $\a  \leq \g _{i_0}$ and $|I_{\a  }|\geq \o $. Let $\la i_k :k\in \N\ra$ be a one-to-one sequence in
$I_{\a  }\setminus \{ i_0\}$ and $I':= \{ i_k :k\in \o \}$. Then
$f\in \Sur (I') \setminus \Sym (I')$, where $f(i_0)=f(i_1)=i_0$ and $f(i_{k+1})=i_k$, for $k\geq 1$.
Now $f^{-1}[\{ i_0\}]=\{ i_0 ,i_1\}$ and,
since $\a  \leq \g _{i_0}\cong \BL _{i_0}'$, by Fact \ref{TC000}(c)
there is a partition $\{ A,B\}$ of $L_{i_0}'$, where $A_{i_1}:= A\cong \a\cong \BL _{i_1}$ and $B\cong \g_{i_0}$.
Thus $A_{i_0}:= B\cup L_{i_0}''\cong \g_{i_0}+n_{i_0}\cong\BL _{i_0}$ and $\{A_{i_0},A_{i_1} \}$ is a partition of $L_{i_0}$.
For $k\geq 1$ we have $f^{-1}[\{ i_k\}]=\{ i_{k+1}\}$ and $\BL _{i_k}\cong \a  \cong \BL _{i_{k+1}}$.
By Fact \ref{TC002} the union $\bigcup _{i\in I'}\BL _i$ is not reversible and, by Fact \ref{TC003}, $\P$ is not reversible.

Let $\g \in \Lim_0$, where $|J_\g \setminus I_\g|\geq \o$ and $\la n_i : i\in J_\g \setminus I_\g\ra$ is not reversible.
Then, by Lemma \ref{TC004}, the union $\bigcup _{i\in J_\g}\BL _i$ is not reversible and, by Fact \ref{TC003}, $\P$ is not reversible.

$(\Rightarrow)$ Suppose that the union $\bigcup_{i\in I}\alpha_i$ is not reversible and that (A) fails, that is,
\begin{equation}\label{EQC001}
\forall i\in I \;\; \forall \a \leq \g _i \;\; |I_\a |<\o.
\end{equation}
By Fact \ref{TB055} we have $K=\{\a \in\Ord: |I_\a |\geq \o \}\neq\emptyset$ and, by (\ref{EQC001})
\begin{equation}\label{EQC002}
\forall i\in I \;\; \forall \a \in K \;\; \g _i <\a ,
\end{equation}
which implies that there is $\g ^*=\max \{ \g _i :i\in I\}\in \Lim _0$
such that for each $i\in \bigcup _{\a \in K}I_\a$ we have $\g _i =\g ^*$, that is, $\a _i=\g ^* +n_i$.

Let $A:=\{ \alpha_i :i\in I \land \alpha _i<\gamma^*\}$ and let $\la \alpha_\xi:\xi<\zeta\ra$ be the increasing enumeration of the set $A$.
Then
$$\textstyle
I=J_{\gamma^*}\cup\,\bigcup_{\xi<\zeta}I_{\alpha_\xi}
$$
is a partition of the set $I$ and, by (\ref{EQC001}),  $|I_{\alpha_\xi}|<\o$, for $\xi<\zeta$.

By Fact \ref{TC002} there is $f\in \Sur (I)\setminus \Sym (I)$ such that for each $j\in I$ the set
$\a_j$ can be partitioned into copies of $\a_i$, for $i\in f^{-1}[\{j\}]$. Since $\a _i \hookrightarrow\alpha_{f(i)}$ we have
\begin{equation}\label{EQC003}
\forall i\in I \;\; \a_i\leq\a_{f(i)}.
\end{equation}
By induction we prove that
\begin{equation}\label{EQC004}
\forall\xi<\zeta\;\;\;f[I_{\alpha_\xi}]=I_{\alpha_\xi}.
\end{equation}
If $j\in I_{\alpha_0}$,
then, by (\ref{EQC003}), for $i\in f^{-1}[\{j\}]$ we have $\alpha_i\leq\alpha_j=\alpha_0$,
which, by the minimality of $\alpha_0$, implies that $\alpha_i=\alpha_0$,
that is $i\in I_{\alpha_0}$.
Thus $f^{-1}[\{j\}]\subset I_{\alpha_0}$, for all $j\in I_{\alpha_0}$,
and hence, $f^{-1}[I_{\alpha_0}]\subset I_{\alpha_0}$.
Since $f$ is onto, we have $I_{\alpha_0}=f[f^{-1}[I_{\alpha_0}]]\subset f[I_{\alpha_0}]$
and, thus, $|I_{\alpha_0}|\leq|f[I_{\alpha_0}]|\leq|I_{\alpha_0}|$.
So, $|I_{\alpha_0}|=|f[I_{\alpha_0}]|$,
which, since the set $I_{\alpha_0}$ is finite and $I_{\alpha_0}\subset f[I_{\alpha_0}]$, implies that $f[I_{\alpha_0}]=I_{\alpha_0}$.

Assuming that $\eta<\zeta$ and $f[I_{\alpha_\xi}]=I_{\alpha_\xi}$, for all $\xi<\eta$,
we prove that  $f[I_{\alpha_\eta}]=I_{\alpha_\eta}$.
If $j\in I_{\alpha_\eta}$ and  $i\in f^{-1}[\{j\}]$, then by (\ref{EQC003})
we have $\alpha_i\leq\alpha_j=\alpha_\eta$.
The inequality $\alpha_i<\alpha_\eta$ would imply that $\alpha_i=\alpha_\xi$, for some $\xi<\eta$,
that is $i\in I_{\alpha_\xi}$, and, by the induction hypothesis, $j=f(i)\in I_{\alpha_\xi}$, which is false.
So, $\alpha_i=\alpha_\eta$ and, hence, $i\in I_{\alpha_\eta}$.
Thus $f^{-1}[\{j\}]\subset I_{\alpha_\eta}$, for all $j\in I_{\alpha_\eta}$,
and hence, $f^{-1}[I_{\alpha_\eta}]\subset I_{\alpha_\eta}$.
Now, as above, we show that $f[I_{\alpha_\eta}]=I_{\alpha_\eta}$ and (\ref{EQC004}) is proved.

By (\ref{EQC004}) and since $|I_{\alpha_\xi}|<\o$, for $\xi<\zeta$,
the restrictions  $f| I_{\alpha_\xi}:I_{\alpha_\xi}\rightarrow I_{\alpha_\xi}$, $\xi<\zeta$, are bijections
and, hence, $f| \bigcup_{\xi<\zeta}I_{\alpha_\xi}:\bigcup_{\xi<\zeta}I_{\alpha_\xi}\rightarrow\bigcup_{\xi<\zeta}I_{\alpha_\xi}$ is a bijection.
By (\ref{EQC003}), we have that $f[J_{\gamma^*}]\subset J_{\gamma^*}$,
and since $f:I\rightarrow I$ is a non-injective surjection,
we conclude that
$f|J_{\gamma^*}:J_{\gamma^*}\rightarrow J_{\gamma^*}$
is a non-injective surjection.
By Fact \ref{TC002} the union $\bigcup_{i\in J_{\gamma^*}}\alpha_i$ is not reversible
and, by Lemma \ref{TC004}, $|I_{\g ^*}|\geq \o$ or the sequence $\langle n_i:i\in J_{\gamma^*}\setminus I_{\g ^*}\rangle$ is not reversible.
But $|I_{\g ^*}|\geq \o$ would imply that $\g ^*\in K$, which is impossible by (\ref{EQC002}). So (B) is true.
\kdok
\noindent
{\bf Proof of Theorem \ref{TC006}.}
It is evident that conditions (I) and (II) are exclusive.

($\Rightarrow$) Let $\P:=\bigcup _{i\in I} \a_i$ be a reversible poset
and suppose that the sequence $\la \a_i :i\in I\ra$ is not finite-to-one, that is, $K=\{\a \in\Ord: |I_\a |\geq \o \}\neq\emptyset$.
By Theorem \ref{TC005} we have $\neg\,$(A) and $\neg\,$(B). Now, as in the proof of Theorem \ref{TC005} we find
$\g ^*\in \Lim _0$ such that for each $i\in \bigcup _{\a \in K}I_\a$ we have $\g _i =\g ^*$, that is, $\a _i=\g ^* +n_i$,
and $\g _i \leq \g ^*$, for all $i\in I$. So $\g ^*=\max \{ \g _i : i\in I\}$.
By $\neg\,$(A) we have $|I_\a|<\o  $, for all $\a \leq \g^*$.
In addition,
$$\textstyle
I=J_{\gamma^*}\cup\,\bigcup_{\xi<\zeta}I_{\alpha_\xi}
$$
is a partition of the set $I$ and  $|I_{\alpha_\xi}|<\o$, for $\xi<\zeta$.
%, which implies that $|I_\a|<\o  $, for all $\a < \g^*$.In addition, by $\neg\,$(A) we have $|I_{\g ^*}|<\o$.
So, if $J_{\gamma^*}=\{ i\in I: \a _i\geq \g^*\}$ would be a finite set,
then the sequence $\la \a_i :i\in I\ra$ would be finite-to-one, which contradicts our assumption.
Thus $|J_{\gamma^*}|\geq \o$, and, since by $\neg$ (A) we have $|I_{\g ^*}|<\o$, it follows that $|J_{\g ^*}\setminus I_{\g ^*}|\geq\o$,
which, together with $\neg\,$(B) implies that  $\la n_i : i\in J_{\g ^*}\setminus I_{\g ^*}\ra$ is a reversible sequence.

($\Leftarrow$) By Fact \ref{TB055}, condition (I) implies that the poset $\P$ is reversible.

Assuming (II) we prove $\neg\,$(A) and $\neg\,$(B). First, if $i\in I$ and $\a \leq \g_i$, then, since $\g_i \leq \g$,
by (II) we have $|I_\a|<\o  $ and $\neg\,$(A) is true.

Second, suppose that (B) is true.
Then there is  $\g '\in \Lim _0$, where $|J_{\g '}\setminus I_{\g '}|\geq\o$ and the sequence $\la n_i : i\in J_{\g '} \setminus I_{\g '}\ra$
is not reversible. By Fact \ref{TC001}, this sequence is not finite-to-one, which means that there is $n\in \N$ such that $|I_{\g ' +n}|\geq \o$.
By (II) we have $\g '<\g$ and, hence, $\a := \g ' +n <\g$. Now, by (II) again $|I_\a|<\o$ and we have a contradiction. Thus we have $\neg\,$(B).
\kdok
\begin{ex}\label{EXC000}\rm
By Theorem \ref{TC006} the poset $\P$ is reversible, since it satisfies (II), where
$$\P =
     \Big(\bigcup _{n\in \N}n \Big)
\cup \Big(\bigcup _{14}\o \Big)
\cup \Big(\bigcup _{\o _1}(\o + 4) \Big)
\cup \Big(\bigcup _{\o _3}(\o + 6) \Big)
\cup \Big(\bigcup _{n\in \o }(\o + 2n+1)\Big).
$$
On the other hand, if
$$
\P_1=\Big(\bigcup _{\o }1\Big) \cup \o \;\;\mbox{ and }\;\;
\P _2=\Big(\bigcup _{\o }(\o +2 )\Big)\cup \Big(\bigcup _{\o }(\o +4)\Big),
$$
then, by Theorem \ref{TC005}, the poset $\P_1$ is not reversible, because it satisfies (A)
and the poset $\P _2$ is not reversible, since it satisfies (B).
\end{ex}
\begin{rem}\label{RC002}\rm
If the ``limit part" $\bigcup_{i\in}\gamma_i$ of the union $\bigcup_{i\in I}\alpha_i$ is reversible,
then $\bigcup_{i\in I}\alpha_i$ is reversible as well.
(By Fact \ref{TC008}(b)  $\langle\gamma_i:i\in I\rangle$ is a finite-to-one sequence, thus $\langle\alpha_i:i\in I\rangle$ is
finite-to-one too and we apply Fact \ref{TB055}).
\end{rem}
\section{Unions of well orders and reversed well orders}\label{S3}
Let $\P =\bigcup _{i\in I}\BL _i$, where $\BL _i\in \CW \cup \CW ^*$, for all $i\in I$, and let
\begin{eqnarray*}
I_{fin} & := & \{ i\in I: |L_i|<\o\},\\
I_w     & := & \{ i\in I: \BL _i\in \CW \land |L_i|\geq \o\} ,\\
I_{w^*} & := & \{ i\in I: \BL _i\in \CW ^*\land |L_i|\geq \o \}.
\end{eqnarray*}
Clearly $\{ I_{w^*} ,I_{fin},I_w\}$ is a partition of $I$. Defining
\begin{eqnarray*}\textstyle
\P_{\CW}& := & \textstyle\bigcup _{i\in I_{fin}\cup I_w}\BL _i \; \mbox{ and }\\[1mm]
\P_{\CW ^*}\! & := & \textstyle\bigcup _{i\in I_{fin}\cup I_{w^*}}\BL _i ,
\end{eqnarray*}
%$$\textstyle
%\P_{\CW}:=\bigcup _{i\in I_{fin}\cup I_w}\BL _i \; \mbox{ and }\; \P_{\CW ^*}:=\bigcup _{i\in I_{fin}\cup I_{w^*}}\BL _i
%$$
we see that $\P_{\CW}$ (resp.\ $\P_{\CW ^*}$) is a disjoint union of well orders (resp.\ reversed well orders)
and the reversibility of these posets is characterized in Theorem \ref{TC006}.
So, if $I_w=\emptyset$ or $I_{w ^*}=\emptyset$, Theorem \ref{TC006} applies.
\begin{te}\label{TC007}
If $\P =\bigcup _{i\in I}\BL _i$, where $\BL _i\in \CW \cup \CW ^*$, for all $i\in I$, and $I_w\neq\emptyset$ and $I_{w ^*}\neq\emptyset$,
then the poset $\P$ is reversible iff the posets $\P_{\CW}$ and $\P_{\CW ^*}$ are reversible.
\end{te}
\dok
($\Rightarrow$) This implication follows from Fact \ref{TC003}.

($\Leftarrow$) Assuming that  $\P_{\CW}$ and $\P_{\CW ^*}$ are reversible posets by Theorem \ref{TC006} we have
the following cases:

{\it Case 1.} At least one of the posets $\P_{\CW}$ and $\P_{\CW ^*}$ satisfies (II).
If $\P_{\CW}$ satisfies (II),
then $\g =\max \{ \g _i : i\in I_w\}\geq \o$ and for $n < \o$ we have $|I_n|<\o $. According to Fact \ref{TC002},
assuming that $f\in \Sur (I)$ and that for each $j\in I$ there exists a partition $\{ A_i: i\in f^{-1}[\{ j\}]\}$ of $L_j$ such that
\begin{equation}\label{EQC005}
\forall i\in f^{-1}[\{ j\}]\;\; \mathbb{L} _i  \cong  \A_i ,
\end{equation}
we show that $f\in \Sym (I)$. Let $\{ n_k :k< \mu \}$ be an increasing enumeration of the set $\{ |L_i|: i\in I \land |L_i|<\o  \}$.
Clearly we have $\mu \leq \o$ and $|I_{n_k}|<\o$, for all $k<\mu$.
Since $f$ is a surjection, an easy induction shows that $f[I_{n_k}]=I_{n_k}$, for all $k<\mu$ and, hence,
$f[I_{fin}]=\bigcup _{k<\mu}f[I_{n_k}]=I_{fin}$. It is evident that $f[I_w]\subset I_w$ and $f[I_{w^*}]\subset I_{w^*}$, which,
since $f$ is onto, implies that $f[I_w]= I_w$ and $f[I_{w^*}]= I_{w^*}$.
Thus $f[I_{fin}\cup I_w]=I_{fin}\cup I_w$, that is,
$$
f_{I_{fin}\cup I_w}:= f\upharpoonright (I_{fin}\cup I_w)\in \Sur (I_{fin}\cup I_w).
$$
In addition, for $j\in I_{fin}\cup I_w$ and $i\in f^{-1}[\{ j\}]$ by (\ref{EQC005}) we have $\BL _i \hookrightarrow \BL _j$, which implies
$i\in I_{fin}\cup I_w$; thus  $f^{-1}[\{ j\}]\subset I_{fin}\cup I_w$ and, hence, $f^{-1}[\{ j\}]= f_{I_{fin}\cup I_w}^{-1}[\{ j\}]$.
Now, since the poset $\P_{\CW}$ is reversible, using Fact \ref{TC002} we obtain
$f\upharpoonright (I_{fin}\cup I_w)\in \Sym (I_{fin}\cup I_w)$, which implies that $f\upharpoonright I_{fin}\in \Sym (I_{fin})$
and $f\upharpoonright I_{w}\in \Sym (I_{w})$.
Similarly we have $f\upharpoonright (I_{fin}\cup I_{w^*})\in \Sym (I_{fin}\cup I_{w^*})$, which implies $f\upharpoonright I_{w^*}\in \Sym (I_{w^*})$.
Thus $f\in \Sym (I)$ indeed.

If $\P_{\CW ^*}$ satisfies (II) we have a similar proof.

{\it Case 2.} The posets $\P_{\CW}$ and $\P_{\CW ^*}$ satisfy (I). Then the sequences of types $\la \otp(\BL _i) :i\in I_{fin}\cup I_w\ra$
and $\la \otp(\BL _i) :i\in I_{fin}\cup I_{w^*}\ra$ are finite-to-one and, hence, the sequence $\la \otp(\BL _i) :i\in I\ra$ is finite-to-one
as well. By Fact \ref{TB055}, the poset $\P$ is reversible.
\kdok
%\section{?}\label{S4}
\begin{rem}\label{RC000}\rm
A topological space $\CX =\la X, {\mathcal O}\ra$ is called reversible iff each continuous bijection $f:\CX \rightarrow \CX$ is a homeomorphism.
The class of reversible spaces contains Euclidean spaces (${\mathbb R}^n$), compact Hausdorff spaces,
and many other classes of relevant topological spaces \cite{RajWil,DoyHoc,Dow}.
If $\P =\la P, \leq \ra$ is a partial order and ${\mathcal O}$ the topology on the set $P$
generated by the base $\{ B_p : p\in P\}$, where $B_p:= \{ q\in p: q\leq p\}$, it is easy to check that endomorphisms of $\P$
are exactly the continuous self mappings of the space $\la P,{\mathcal O}\ra$. Thus,
the poset $\P$ is reversible iff $\la P,{\mathcal O}\ra$ is a reversible topological space.
So, the results of this article can be interpreted as characterizations of reversibility in the class
of topological spaces of the form $\la P,{\mathcal O}\ra$, where $\P =\bigcup _{i\in I}\BL _i$ and $\BL _i\in \CW \cup \CW ^*$, for $i\in I$.
%Clearly these spaces are disconnected $T_0$.
\end{rem}
\begin{rem}\label{RC001}\rm
Clearly, the posets of the form $\P =\bigcup _{i\in I}\BL _i$, where $\BL _i\in \CW $, for $i\in I$, are, in fact, trees.
Thus, Theorem \ref{TC006} can be understood as a characterization of reversibility in the class of such trees
and using it we can construct reversible trees  of arbitrary size and height having all levels infinite.
We recall that, by \cite{Kuk},  well founded posets with finite levels are reversible.
\end{rem}
The posets investigated in this paper are disjoint unions of scattered linear orders and it is natural to ask what is going on
when some of $\BL _i$'s are not scattered. The following example is related to that question.
\begin{ex}\label{EXC001}\rm
If $\P =\bigcup _{i\in \o}\BL _i$, where $\BL _i$, $i\in \o$, are pairwise disjoint countable linear orders,
$\BL _0=\Q$ and $\BL _i$ is a CSB linear order, for $i>0$, we prove that
$$
\P \mbox{ is reversible } \Leftrightarrow \la \otp (\BL _i): i\in \o\ra \mbox{ is a finite-to-one sequence}.
$$
Since countable CSB linear orders are scattered, for $i>0$ we have $\BL_i \in \Scatt$ and $\otp (\BL _i)= [\BL _i]_{\rightleftarrows}$.
So if $\la \otp (\BL _i): i\in \o\ra$ is a finite-to-one sequence, the sequence $\la [\BL _i]_{\rightleftarrows}: i\in \o\ra$ is finite-to-one as well
and, by Fact \ref{TB055}, the poset $\P $ is reversible.

Conversely, if the sequence $\la \otp (\BL _i): i\in \o\ra$ is not finite-to-one  and $|I_\t|=\o$, for some order-type $\t$,
then $\t$ is a scattered type. Let $i_0=0$ and let $\la i_k : k\in \N \ra$ be a one-to-one sequence in $I_\t$.
Let $g : \BL _{i_1}\hookrightarrow \Q =\BL _{i_0}$. It is easy to see that $Q \setminus g[L _{i_1}]$ is a dense linear order without end-points
and, by Cantor's theorem, $\Q \cong Q \setminus g[L _{i_1}]$.
Now, defining $f\in \Sur ( \{ i_k : k\in \o \})\setminus \Sym ( \{ i_k : k\in \o \})$ by
$f(i_0)=f(i_1)=i_0$ and $f(i_k)=i_{k-1}$, for $k\geq 2$, using Fact \ref{TC002} we conclude that the poset $\bigcup _{k\in \o}\BL _{i_k}$ is not reversible,
which by Fact \ref{TC003} implies that $\P $ is not reversible too. The claim is proved

Of course the statement given above will hold in a slightly more general situation, when we have finitely many countable non-scattered
chains, one of them contains a convex copy of $\Q$ and all but finitely many scattered components are countable CSB linear orders.

We note that, in particular, the statement holds if the countable scattered $\BL _i$'s are from $\CW \cup \CW^*$.
Then, regarding Theorem \ref{TC006}, the ``finite-to-one" condition has no alternative (case (II) cannot appear).
\end{ex}
\section{Reversible sequences of several things}
Here we discuss several concepts of a ``reversible sequence" related to this article.
First, condition (\ref{EQC007}) defines a {\it reversible sequence of structures};
(\ref{EQB032}) defines a {\it reversible sequence of natural numbers} and, in \cite{KuMoEqu},
more generally, a sequence of non-zero cardinals $\langle\kappa_i:i\in I\rangle$ is called
a {\it reversible sequence of cardinals} iff
\begin{equation}\label{EQC008}\textstyle
\neg\exists f\in\Sur(I)\setminus\Sym(I)\;\;\;\forall j\in I\;\sum_{i\in f^{-1}[\{j\}]}\kappa_i=\kappa_j.
\end{equation}
In addition, defining a sequence of non-zero ordinals $\langle\alpha_i:i\in I\rangle$ to be a {\it reversible  sequence of ordinals} iff
it satisfies (I) or (II) of Theorem \ref{TC006}, we have
\begin{prop}\label{TC009}
For each sequence of non-zero cardinals $\bar{\kappa}=\langle\kappa_i:i\in I\rangle$ the following conditions are equivalent:

(a) The poset $\bigcup_{i\in I}\k _i$ is a reversible structure,

(b) $\bar{\kappa}$ is a reversible sequence of structures,

(c) $\bar{\kappa}$ is a reversible sequence of ordinals,

(d) $\bar{\kappa}$ is a reversible sequence of cardinals,

(e) $\bar{\kappa}$ is a finite-to-one sequence or a reversible sequence of natural numbers.
\end{prop}
\dok
The equivalence (a) $\Leftrightarrow$ (b) follows  from Fact \ref{TC002},  (a) $\Leftrightarrow$ (c) is Theorem \ref{TC006}
and (d) $\Leftrightarrow$ (e) is Theorem 1.1 of \cite{KuMoEqu}.

(d) $\Rightarrow$ (b) We prove the contrapositive: if $f\in\Sur(I)\setminus\Sym(I)$
and for each $j\in I$ we have $\bigcup_{i\in f^{-1}[\{j\}]}\kappa_i\preccurlyeq_c\kappa_j$,
then, clearly, $\sum_{i\in f^{-1}[\{j\}]}\kappa_i=\kappa_j$, for $j\in J$.

(c) $\Rightarrow$ (e) If $\bar{\kappa}$ satisfies (I), then (e) is true. If $\bar{\kappa}$ satisfies (II),
then there is $\gamma=\max\{\gamma_i:i\in I\}$
and the sequence $\la n_i :i\in J_\g \setminus I_\g\ra$ is not finite-to-one,
which implies that $|J_\gamma\setminus I_\gamma|\geq \o$.
For $i\in J_\gamma\setminus I_\gamma$ we have $\g _i=\g \neq \k _i =\g _i +n_i$,
which implies that $\gamma=0$.
Thus $\bar{\kappa}\in{}^{I}\mathbb{N}$ and $I_\gamma=\emptyset$, which implies that $J_\gamma\setminus I_\gamma=I$, and, by (II),
$\bar{\kappa}=\la n_i:i\in I\ra$ is a reversible sequence of natural numbers.
\kdok
\begin{rem}\rm
By \cite{KuMoEqu}, assuming that $\la \X_i \!:\!i\in I\ra$ is a sequence of pairwise disjoint, connected and reversible binary structures
which is {\it rich for monomorphisms} (RFM), which means that for all $i,j\in I$ and each $A\in[X_j]^{|X_i|}$ there is a monomorphism $g:\mathbb{X}_i\rightarrow\mathbb{X}_j$ such that $g[X_i]=A$, we have:
\begin{equation}\label{EQC009}\textstyle
\bigcup_{i\in I}\mathbb{X}_i \mbox{ is reversible}
\Leftrightarrow \langle |X_i|:i\in I\rangle \mbox{ is a reversible sequence of cardinals}.
\end{equation}
In particular, this characterizes reversible equivalence relations and disjoint unions of cardinals $\leq \o$.
We remark that the equivalence (a) $\Leftrightarrow$ (d) of Proposition \ref{TC009} shows that the characterization (\ref{EQC009})
holds in a class of (sequences of) structures which is larger than the class RFM
(for example, $\la \o ,\o _1,\o _2,\o _3, \dots\ra \not\in \rfm$).
\end{rem}
\begin{prop}\label{TC010}
For each sequence of non-zero ordinals $\bar{\alpha}=\langle\alpha_i:i\in I\rangle$ the following conditions are equivalent:

(a) The poset $\bigcup_{i\in I}\a _i$ is a reversible structure,

(b) $\bar{\alpha}$ is a reversible sequence of structures,

(c) $\bar{\alpha}$ is a reversible sequence of ordinals.\\
In addition, any of the conditions listed above is implied by the following:

(d) $\langle |\alpha_i|:i\in I\rangle$ is a reversible sequence of cardinals.
\end{prop}
\dok
The equivalence (a) $\Leftrightarrow$ (b) follows  from Fact \ref{TC002},  (a) $\Leftrightarrow$ (c) is Theorem \ref{TC006}.

If (d) is true, then, by Proposition \ref{TC009}, either $\langle |\alpha_i|:i\in I\rangle$ is a finite-to-one sequence, thus
$\langle \alpha_i:i\in I\rangle$ is finite-to-one, or  $\langle |\alpha_i|:i\in I\rangle=\langle \alpha_i:i\in I\rangle\in {}^{I}\N$
is a reversible sequence of natural numbers. So, by Theorem \ref{TC006}, (a) is true.
\hfill $\Box$
\begin{ex}\rm
(c) $\not\Rightarrow$ (d) in Proposition \ref{TC010}:
$\langle\omega+n:n\in\omega\rangle$ is a reversible sequence of ordinals,
but $\langle\omega:n\in\omega\rangle$ is not a reversible sequence of cardinals.
\end{ex}
\paragraph{Acknowledgments}
This research was supported by the Ministry of Education and Science of the Republic of Serbia (Project 174006).


\footnotesize
\begin{thebibliography}{99}
\bibitem{DoyHoc}
	  P.\ H.\ Doyle, J.\ G.\ Hocking,
	  Bijectively related spaces,
	  I. Manifolds. Pac. J. Math. 111 (1984) 23--33.
\bibitem{Dow}
      A.\ Dow, R.\ Hern\'{a}ndez-Guti\'{e}rrez,
      Reversible filters,
      Topology Appl.\ 225 (2017) 34--45.
\bibitem{Kuk}
	  M.\ Kukie\strokel a,
	  Reversible and bijectively related posets,
	  Order 26 (2009) 119--124.
\bibitem{Kuk1}
      M.\ Kukie\strokel a,
      Characterization of hereditarily reversible posets,
      Math.\ Slovaca 66,3 (2016) 539--544.
\bibitem{KRet}
      M.\ S.\ Kurili\'c,
      Retractions of reversible structures, J.\ Symbolic Logic, (in print).
\bibitem{KDef}
   	  M.\ S.\ Kurili\'c,
	  Reversibility of definable relations,
	  (to appear)
\bibitem{KuMoSim}
      M.\ S.\ Kurili\'c, N.\ Mora\v ca,
      Condensational equivalence, equimorphism, elementary equivalence and similar similarities,
      Ann.\ Pure Appl.\ Logic 168,6 (2017) 1210--1223.
\bibitem{KuMoExtr}
      M.\ S.\ Kurili\'c, N.\ Mora\v ca,
	  Reversibility of extreme relational structures,
	  (to appear)
\bibitem{KuMoEqu}
	  M.\ S.\ Kurili\'c, N.\ Mora\v ca,
	  Reversible sequences of cardinals, reversible equivalence relations, and similar structures,
	  (to appear) https://arxiv.org/abs/1709.09492
\bibitem{KuMoDiscI}
	  M.\ S.\ Kurili\'c, N.\ Mora\v ca,
	  Reversibility of disconnected structures,
	  (to appear)
\bibitem{Laf}
      C.\ Laflamme, M.\ Pouzet, R.\ Woodrow,
      Equimorphy: the case of chains,
      Arch.\ Math.\ Logic 56, 7--8 (2017) 811--829.
\bibitem{Lav1}
      R.\ Laver,
      An order type decomposition theorem
      Ann.\ of Math.\ 98,1 (1973) 96--119.
\bibitem{RajWil}
	  M.\ Rajagopalan, A.\ Wilansky,
	  Reversible topological spaces,
	  J.\ Aust.\ Math.\ Soc.\ 61 (1966) 129--138.
\bibitem{Rosen}
      J.\ G.\ Rosenstein,
      Linear orderings,
      Pure and Applied Mathematics, 98, Academic Press, Inc.,
      Harcourt Brace Jovanovich Publishers, New York-London, 1982.
\end{thebibliography}
\end{document}